# Positive Dual Certificates for Real Separated Subset Sums: Chamber Optimality of the Conway–Guy and Bohman Sets, and Exact Values up to Eight

Esther Notik

Independent researcher, Petah Tikva, Israel — esthnt@gmail.com

## Abstract

Let $f_R(n)$ be the least possible largest element of a set of n positive reals whose $2^n$ subset sums are pairwise at least 1 apart, and let $f_Z(n) \geq f_R(n)$ be the corresponding integer quantity of the Erdős distinct subset sums problem. Fixing the signs of all signed relations $\varepsilon \cdot x$, $\varepsilon \in \{-1,0,1\}^n$, decomposes the real problem into polyhedral *chambers*, on each of which it is a linear program. We prove that for every n the Conway–Guy set is the unique optimum of its chamber. The proof is an explicit positive dual certificate consisting of n signed relations that are tight at the Conway–Guy set and unit lower triangular in gap coordinates; the relations are described by blocks of consecutive coordinates together with the positions $n - k(k+1)/2$, and tightness reduces to an identity of the Conway–Guy sequence proved by induction. We then show, by a branch-and-bound over chambers in which every leaf carries an exact rational certificate, that $f_R(n) = f_Z(n)$ for all $n \leq 8$; in particular $f_R(7) = 44$ and $f_R(8) = 84$. We then prove the same for both of Bohman's families, which first improve on Conway–Guy at n = 12 (1159 < 1164): every Bohman set is the unique optimum of its chamber, with a block-triangular certificate given by explicit formulas. Finally we exhibit chambers whose optimal vertex is not an integer point.

## 1. Introduction

A set $\{a_1, \ldots, a_n\}$ of positive integers has *distinct subset sums* if the $2^n$ sums $\Sigma_{i \in S}\, a_i$, $S \subseteq \{1, \ldots, n\}$, are pairwise different. Let $f_Z(n)$ denote the least possible value of max $a_i$. Erdős conjectured that $f_Z(n) \geq c \cdot 2^n$ for an absolute constant c > 0 [E].

A disproof was reported in 2026: for every $\varepsilon > 0$ there are arbitrarily large n with $f_Z(n) \leq \varepsilon \cdot 2^n$. The argument is ineffective (it gives no bound on how large n must be); it was found by an AI system in the FrontierMath Erdős benchmark of Adamczewski and Bloom and is machine-checked in Lean 4 [AB26, Ad26]. Since $f_R(n) \leq f_Z(n)$, the same holds for the real relaxation. Explicit sum-distinct sets below Bohman's constant 0.22002 have since been reported by Toledo [To26], who made the transfer step of the proof effective, building on Bloom's lattice interpretation; his smallest example has n = 1,701 elements and $N/2^n = 0.217967$. The results below are finite and are not affected by this development.

The best general lower bounds are of order $2^n/\sqrt{n}$ (Erdős–Moser; Elkies; Dubroff, Fox and Xu [DFX21], who proved $f_Z(n) \geq C(n, \lfloor n/2 \rfloor)$; Steinerberger [St22]).

Classical explicit constructions include those of Conway and Guy [CG68], proved to have distinct subset sums for every n by Bohman [B96, Theorem 1], and the families of Lunnon [L88] and Bohman [B96, Theorem 3; B98]. Bohman's sets first improve on the Conway–Guy sets at n = 12, where they give 1159 instead of 1164 (Section 4.1).

Exact values are known for $n \leq 10$, where the Conway–Guy sets are optimal: Lunnon determined $f_Z(n)$ for $n \leq 8$ by exhaustive search [L88], and Borwein and Mossinghoff treated n = 9 [BM03] (the OEIS entry A276661 credits $f_Z(9) = 161$ to J. P. Grossman [OEIS]). P. W. Dyson found $f_Z(10) = 309$ in 2025 by an exhaustive integer search, with the Conway–Guy set as the only optimal set [Dy25]; this computation has not been refereed and comes without a certificate.

The lower-bound arguments use only that distinct subset sums are separated by at least 1, and so they apply to the following real relaxation, which is the subject of this paper.

**Definition 1.1.** *For $n \geq 1$ let $f_R(n)$ be the infimum of $\max a_i$ over all n-tuples of positive reals such that $|\Sigma_{i\in S} a_i - \Sigma_{i\in T} a_i| \geq 1$ for all subsets $S \neq T$.*

Every integer set with distinct subset sums is admissible, so $f_R(n) \leq f_Z(n)$. The basic question is whether the inequality is ever strict, and if so for which first n.

Costa and Della Fiore [CDF26, Theorem 1.1 and Proposition 2.2] recently proved $f_Z(n) \geq H_n$ by ordering all $2^n$ subsets by their sums and applying Harper's formula for the bandwidth of the hypercube [H66; see also WWD09]; this improves the bound $C(n, \lfloor n/2 \rfloor)$ of Dubroff, Fox and Xu for every $n \geq 3$. Their argument is stated for integers, but it uses only that distinct subset sums are at least 1 apart, so it applies to the real relaxation. We record this observation.

**Proposition 1.2.** *$f_R(n) \geq H_n := \Sigma_{j=0}^{n-1} C(j, \lfloor j/2 \rfloor)$. In particular $f_R(8) \geq 78$.*

*Proof* (following [CDF26, Proposition 2.2]). Let a be admissible and number the vertices of the n-cube $Q_n$ by the increasing order of their subset sums w(S). By Harper's theorem the bandwidth of $Q_n$ equals $H_n$, so some edge $\{S, S \cup \{i\}\}$ joins positions $d \geq H_n$ apart. The $d + 1$ sums from w(S) to $w(S \cup \{i\}) = w(S) + a_i$, inclusive, are pairwise at least 1 apart, so $a_i \geq d \geq H_n$. In the integer case [CDF26] counts the integers in the open interval $(w(S), w(S) + a_i)$; for reals the 1-separation replaces this count. ■

The values are $H_6 = 23$, $H_7 = 43$, $H_8 = 78$, $H_9 = 148$. Thus $78 \leq f_R(8) \leq f_Z(8) = 84$.

### *Main results*

- **Theorem A (Section 3).** For every n, the Conway–Guy n-set is the unique minimiser of $\max x_i$ over its sign chamber. The certificate is explicit and uniform in n.
- **Theorem B (Section 5).** $f_R(n) = f_Z(n)$ for $1 \leq n \leq 8$. In particular $f_R(7) = 44$ and $f_R(8) = 84$. The proof is computer-assisted, with exact rational certificates checked by a program independent of the search code; the case n = 7 was also checked by a second, separately written checker, made with the use of Astra (OpenAI).
- **Theorem C (Section 4).** For both of Bohman's families and every $k \geq 2$, every set of the family is the unique minimiser over its sign chamber, again with an explicit positive certificate. These sets first improve on the Conway–Guy sets at n = 12 (Family I: 1159 < 1164) and n = 13 (Family II: 2262 < 2284). The proof extends the induction of Theorem A by a propagation lemma for general rule sequences, together with a finite seed and a top block given by explicit formulas in k.
- **Computational results (Section 6).** Some chambers have optimal vertices that are not integer points.

The method is linear programming duality. Its point is that every lower bound is witnessed by a short positive combination of signed relations, which can be checked by hand or by exact arithmetic.

## 2. Chambers, certificates and gap coordinates

A *signed relation* is a linear form $R_\varepsilon(x) = \varepsilon\cdot x$ with $\varepsilon \in \{-1,0,1\}^n \setminus \{0\}$. A vector x is admissible (1-separated) if and only if $|R_\varepsilon(x)| \geq 1$ for every ε, because a difference of two subset sums is exactly a signed relation. Without loss of generality $0 < x_1 < \cdots < x_n$, so that $\max x_i = x_n$; the relations $x_1$ and $x_i - x_{i-1}$ are then positive.

**Definition 2.1.** *Let a be admissible. Its sign chamber is the polyhedron $C(a) = \{x \in \mathbb{R}^n : \sigma_\varepsilon(a)\cdot R_\varepsilon(x) \geq 1$ for all $\varepsilon\}$, where $\sigma_\varepsilon(a) = \operatorname{sign} R_\varepsilon(a)$.*

The admissible set is the union of the chambers, and on each chamber the minimisation of $x_n$ is a linear program with integer data. There are finitely many chambers, each LP is bounded below, and its feasible region contains no line (it lies in $x_1 \geq 1$, $x_i - x_{i-1} \geq 1$). An optimum is attained at a vertex. Vertices are rational, so $f_R(n)$ is attained at a rational point, and after clearing denominators

$$f_R(n) = \min \{ \max_i a_i / \mathrm{sep}(a) : a \in Z_{>0}^n, \mathrm{sep}(a) > 0 \}, \quad \mathrm{sep}(a) = \min_{\varepsilon \neq 0} |\varepsilon \cdot a|.$$

**Definition 2.2.** *A positive certificate for the bound $x_n \geq c$ on a polyhedron $\{x : R_j(x) \geq 1, j = 1, \dots, N\}$ is a vector $y \in \mathbb{R}^N_{\geq 0}$ with $\Sigma_j y_j R_j(x) = x_n$ identically and $\Sigma_j y_j = c$.*

**Gap coordinates.** Put $d_1 = x_1$ and $d_i = x_i - x_{i-1}$. If $E_i = \Sigma_{k \geq i} \varepsilon_k$ are the suffix sums of ε, then $R_\varepsilon(x) = \Sigma_i E_i d_i$, and $x_n = d_1 + \cdots + d_n$. This reduction is due to Bohman [B96, Lemma 1]: in his terms, $(E_n, E_{n-1}, \dots, E_1)$ is a *smooth* vector ($|v(1)| \leq 1$ and $|v(i) - v(i+1)| \leq 1$), and a set has distinct subset sums if and only if no nonzero smooth vector is orthogonal to its difference vector. Bohman uses it to prove that no relation vanishes; we use the same coordinates to describe relations that equal exactly 1 and to build positive certificates.

## 3. The Conway–Guy chambers

### 3.1 Definitions

The Conway–Guy sequence [CG68] is, in the form used by Bohman [B96, §2], $u_0 = 0$, $u_1 = 1$, $u_{k+1} = 2u_k - u_{k-r(k)}$, with $r(k) = \mathrm{round}(\sqrt{(2k)})$; it begins 0, 1, 2, 4, 7, 13, 24, 44, 84, 161, 309, 594, 1164, … The Conway–Guy n-set is $a_i = u_n - u_{n-i}$ ($1 \leq i \leq n$); we put $a_0 := 0$. Let $T_k = k(k+1)/2$.

For $m \geq 1$ let $p(m)$ be the largest $p \geq 2$ with $p + T_{p-2} < m$, and $p(m) = 1$ if no such p exists. Equivalently, $p(1) = p(2) = 1$ and, for $p \geq 2$, $p(m) = p$ if and only if $(p^2 - p + 4)/2 \leq m \leq (p^2 + p + 2)/2$.

**Definition 3.1.** *For $n \geq 1$ and $1 \leq j \leq n$ put $m = n - j + 1$, $p = p(m)$, and (with $x_0 := 0$)*

$$R^{(j)}(x) = (x_j + x_{j+1} + \cdots + x_{j+p-1}) - x_{j-1} - \Sigma_{k=0}^{p-2} x_{n-T_k}.$$

The condition $p + T_{p-2} < m$ is equivalent to $j + p - 1 < n - T_{p-2}$. Hence the block of p plus signs lies strictly below the positions $n, n - 1, n - 3, \dots, n - T_{p-2}$ carrying minus signs, the coefficients are in $\{-1, 0, 1\}$, and $R^{(j)}$ is a signed relation.

**Theorem 3.2 (Theorem A).** *For every $n \geq 1$ there are integers $y_1, \dots, y_n \geq 1$ such that (i) $\Sigma_j y_j R^{(j)}(x) = x_n$ for all $x \in \mathbb{R}^n$, and (ii) $R^{(j)}(a) = 1$ for every j. Consequently, for every x with $R^{(j)}(x) \geq 1$ ($1 \leq j \leq n$),*

$$x_n - u_n = \Sigma_j y_j (R^{(j)}(x) - 1) \geq 0,$$

*with equality if and only if $x = a$. In particular the Conway–Guy set is the unique minimiser of $x_n$ over its chamber $C(a)$, and the minimum is $u_n$.*

### 3.2 Proof

**Step 1: triangularity.** Let $E^{(j)}_i$ be the suffix sums of the coefficient vector of $R^{(j)}$. For $i < j$ the term $-x_{j-1}$ cancels one plus sign, so $E^{(j)}_i = 0$; and $E^{(j)}_j = p - (p - 1) = 1$. Inside the block, $E^{(j)}_i = j + 1 - i \leq 0$ for $i > j$, and above the block $E^{(j)}_i$ is minus the number of top minus positions $\geq i$. Thus $E^{(j)}_i \leq 0$ for all $i > j$. The identity $\Sigma_j y_j R^{(j)} = x_n$ is equivalent to $\Sigma_j y_j E^{(j)}_i = 1$ for all i, a unit lower triangular system with solution

$$y_1 = 1, \quad y_i = 1 + \Sigma_{j<i} y_j |E^{(j)}_i|.$$

All $y_i$ are integers and $y_i \geq 1$.

**Step 2: tightness.** Each $R^{(j)}$ has as many plus as minus signs (counting $x_0$ if $j = 1$). Substituting $a_i = u_n - u_{n-i}$, the terms $u_n$ cancel and

$$R^{(j)}(a) = F(m) := u_m - \Sigma_{t=1}^{p} u_{m-t} + \Sigma_{k=0}^{p-2} u_{T_k}, \quad p = p(m).$$

**Lemma 3.3.** *(a) If $(p^2 - p + 4)/2 \le N \le (p^2 + p)/2$, then $r(N) = p$. (b) If $N = (p^2 + p + 2)/2$, then $r(N) = p + 1$.*

*Proof.* (a) $(p - ½)^2 < p^2 - p + 4 \le 2N \le p^2 + p < (p + ½)^2$. (b) $(p + ½)^2 < p^2 + p + 2 < (p + 3/2)^2$. Since $\sqrt{(2N)}$ is never a half-integer, rounding is unambiguous. ■ Lemma 3.3 is a precise form of Bohman's observation [B96, §2] that, for $b(n) = r(n - 1)$, each positive integer p occurs in the sequence b exactly p times: $r(N) = p$ exactly for $(p^2 - p + 2)/2 \le N \le (p^2 + p)/2$.

**Lemma 3.4.** *$F(m) = 1$ for every $m \ge 1$.*

*Proof.* $F(1) = u_1 - u_0 = 1$ and $F(2) = u_2 - u_1 = 1$. Let $m \ge 3$. If $p(m) = p(m - 1) = p$, then $F(m) - F(m - 1) = u_m - 2u_{m-1} + u_{m-1-p} = u_{m-1-p} - u_{m-1-r(m-1)}$, which is 0 by Lemma 3.3(a). If $p(m) = p + 1$ and $p(m - 1) = p$, then $m - 1 = (p^2 + p + 2)/2 = T_{p-1} + p + 1$ and $F(m) - F(m - 1) = u_m - 2u_{m-1} + u_{T_{p-1}} = u_{T_{p-1}} - u_{m-1-r(m-1)}$, which is 0 by Lemma 3.3(b). ■

**Step 3: conclusion.** By (i) and $\Sigma_j y_j = \Sigma_j y_j R^{(j)}(a) = a_n = u_n$, we get $x_n - u_n = \Sigma_j y_j (R^{(j)}(x) - 1)$. On C(a) every $R^{(j)}(x) \ge 1$, so $x_n \ge u_n$. If $x_n = u_n$, positivity of the $y_j$ forces $R^{(j)}(x) = 1$ for all j; since the system is unit triangular in gap coordinates, $x = a$. ■

**Remarks.** (1) Since $r(N) = p$ exactly for $T_{p-1} + 1 \le N \le T_p$, the rounding rule changes at $N = T_p + 1$, not at $N = T_p$. At the corresponding transition $m - 1 = T_p + 1$ of Lemma 3.4, the new tail index is $m - 1 - r(m - 1) = T_p - p = T_{p-1}$. This shift is the reason why the minus positions in Definition 3.1 are $n - T_k$.

(2) For $n \le 15$ an exhaustive search shows that $R^{(j)}$ is the only tight relation with the triangular gap pattern of Step 1, and for $5 \le n \le 11$ the certificate coincides with the one selected by the LP that minimises the weighted sum $\Sigma_j w_j y_j$ over all certificates supported on tight relations, where $w_j$ is the number of nonzero coefficients of the j-th relation.

(3) The identities of Theorem 3.2 were also checked in exact integer arithmetic for $2 \le n \le 150$.

(4) For n = 7 the certificate is

$$16(a_7 - a_6) + 13(a_6 - a_5) + 7(a_5 + a_6 - a_4 - a_7) + 4(a_4 + a_5 - a_3 - a_7)$$
$$+ 2(a_3 + a_4 + a_5 - a_2 - a_6 - a_7) + (a_2 + a_3 + a_4 - a_1 - a_6 - a_7) + (a_1 + a_2 + a_3 - a_6 - a_7) = a_7.$$

| **n** | **$u_n$** | **multipliers $y_1, \ldots, y_n$** |
|---|---|---|
| 8 | 84 | 1, 1, 2, 4, 8, 13, 25, 30 |
| 9 | 161 | 1, 1, 2, 4, 8, 15, 25, 48, 57 |
| 10 | 309 | 1, 1, 2, 4, 8, 15, 29, 48, 92, 109 |
| 11 | 594 | 1, 1, 2, 4, 8, 15, 29, 56, 92, 177, 209 |
| 12 | 1164 | 1, 1, 2, 4, 8, 16, 29, 57, 110, 180, 347, 409 |

Theorem 3.2 is local. It cannot be global for all n, because already at n = 12 a Bohman set has a smaller maximum ($1159 < 1164 = u_{12}$; Section 4).

## 4. Bohman chambers

### *4.1 The two families*

Bohman [B96, §4, Theorem 3; B98, Theorems 2.1 and 2.2], extending sequences found by Lunnon [L88], defines two families of difference sequences. By [B98, §2], the four sequences of Lunnon that improve on the Conway–Guy sequence are Families I and II with k = 2 and k = 3, so Theorems 4.3 and 4.4 below cover them. Let $k \ge 1$.

- **Family I.** $d(k-j) = 2\cdot 4^{j-1}$ $(1 \le j \le k-1)$, $d(k) = 1$, $d(k+j) = 4^{j-1}$ $(1 \le j \le k)$; rule $b(2k+1) = b(2k+2) = k+1$, $b(2k+3) = k+2$, and $b(i) = \mathrm{round}(\sqrt{(2(i+2-2k))})$ for $i \ge 2k+4$.
- **Family II.** $d(k+1-j) = 4^{j-1}$, $d(k+1) = 1$, $d(k+1+j) = 2\cdot 4^{j-1}$ $(1 \le j \le k)$; rule $b(2k+2) = k+1$, $b(2k+3) = b(2k+4) = k+2$, and $b(i) = \mathrm{round}(\sqrt{(2(i+1-2k))})$ for $i \ge 2k+5$.

In both families $d(i) = \Sigma_{j=i-b(i)}^{i-1} d(j)$ for every i beyond the initial values, and the sets are $S_m = \{\Sigma_{j=i}^{m} d(j) : 1 \le i \le m\}$; Bohman proved that they have distinct subset sums. Family II with k = 1 is the Conway–Guy sequence. As in Section 3 put $U(t) = d(1) + \cdots + d(t)$ and $U(0) = 0$; then the sorted set is $a_t = U(m) - U(m-t)$, $1 \le t \le m$, and we put $a_0 := 0$.

*Remark.* The example $d_{3,2}$ = (16, 4, 1, 1, 2, 8, 32, 43, 86, 171, 200, 400, …) printed in [B96, §4] does not follow from the definition, which gives (16, 4, 1, 1, 2, 8, 32, 43, 86, 171, 300, 600, …). We use the definition; distinct subset sums were also checked directly for $k \le 5$ and $m \le 22$.

| m | max $S_m$ (Family I, k = 2) | $u_m$ (Conway–Guy) |
|---|---|---|
| 12 | 1159 | 1164 |
| 13 | 2271 | 2284 |
| 14 | 4456 | 4484 |
| 15 | 8748 | 8807 |
| 16 | 17182 | 17305 |

For $n \le 11$ no set of either family has a smaller maximum than the Conway–Guy set. Family I first improves on it at n = 12 (k = 2, maximum 1159 < 1164 = $u_12$), and Family II at n = 13 (k = 2, maximum 2262 < 2284 = $u_13$). Family I with k = 2 stays below $u_n$ for all $12 \le n \le 100$ (checked).

### *4.2 Relations in U-form and a propagation lemma*

A balanced signed relation (coefficient sum 0, the term $a_0$ included when needed) can be written as $V = \Sigma_t c_t U(t)$: under $a_t = U(m) - U(m-t)$, the coefficient of $a_{m-t}$ is $-c_t$. Its value at $S_m$ is $\Sigma_t c_t U(t)$, which does not depend on m. We write $\Sigma[a..b] = U(a) + U(a+1) + \cdots + U(b)$, with $\Sigma[a..b] = 0$ if $a > b$.

Let $i_0$ be an index, let τ be a U-form with indices smaller than $i_0 - b(i_0)$, and put $V_{i0} = U(i_0) - \Sigma[i_0 - b(i_0)..i_0 - 1] + \tau$. For $m' > i_0$ define

$$V_{m'} = U(m') - \Sigma[m' - b(m')..m' - 1] + \tau + \Sigma_{i0 < i \le m',\, b(i) = b(i-1)+1} U(i - 1 - b(i)).$$

**Lemma 4.1 (propagation).** *Suppose that $d(i) = \Sigma_{j=i-b(i)}^{i-1} d(j)$ and $b(i) - b(i-1) \in \{0, 1\}$ for every $i > i_0$. Then $V_{m'}$ and $V_{i0}$ have the same value for every $m' \ge i_0$.*

*Proof.* The recursion gives $U(i) - 2U(i-1) = -U(i-1-b(i))$. If $b(m') = b(m'-1) = p$, then $V_{m'} - V_{m'-1} = U(m') - 2U(m'-1) + U(m'-1-p) = 0$. If $b(m') = b(m'-1) + 1$, the block gains one term and the tail gains $U(m'-1-b(m'))$, so $V_{m'} - V_{m'-1} = U(m') - 2U(m'-1) + U(m'-1-b(m')) = 0$. ■

Theorem 3.2 is the special case $b(i) = \mathrm{round}(\sqrt{(2(i-1))})$, $i_0 = 2$, $\tau = 0$; there Lemma 3.3 identifies the block length with b.

**Lemma 4.2 (validity).** *In addition to the hypotheses of Lemma 4.1, suppose: (i) no two consecutive indices $i > i_0$ are increases of b; (ii) the coefficients of τ are ±1, their sum is $b(i_0) - 1$, and every index of τ is smaller than the first new index $i - 1 - b(i)$; (iii) for every T, the number of coefficients +1 of τ at indices ≤ T is at least the number of coefficients −1 at indices ≤ T. Then for $i_0 \le m' \le m$, $V_{m'}$ is a signed relation for $S_m$, and with $j = m - m' + 1$ its suffix-sum vector satisfies $E_i = 0$ $(i < j)$, $E_j = 1$, $E_i \le 0$ $(i > j)$.*

*Proof.* Between two increases of b the quantity $i - b(i)$ grows by one at each step without increase and stays constant at an increase; by (i) there is at least one step without increase between two increases.

Hence the new indices $i - 1 - b(i)$ are strictly increasing, and they are smaller than $m' - b(m')$, the lowest index of the block. With (ii), all indices of $V_{m'}$ are distinct, every coefficient is $\pm 1$, and the coefficient sum is $1 - b(m') + (b(m') - 1) = 0$. In a-coordinates the head $U(m')$ gives $-1$ at position $j - 1$ (nothing if $j = 1$), the block gives $+1$ at positions $j, \ldots, j + b(m') - 1$, and a tail coefficient $c$ at index $t$ gives $-c$ at position $m - t > j + b(m') - 1$. Therefore $E_i = 0$ for $i < j$ (balance), $E_j = b(m') - (b(m') - 1) = 1$, $E_i = j + 1 - i \le 0$ inside the block, and above the block $E_i$ is minus the sum of the tail coefficients at indices $\le m - i$, which is $\le 0$ by (iii); new tail terms are $+1$ and keep (iii). ■

### *4.3 Chamber optimality of both families*

**Theorem 4.3 (Family I).** *Let $k \ge 2$ and $m \ge 2k$. The set $S_m$ of Family I is the unique minimiser of $x_m$ over its sign chamber. A certificate consists of the m relations listed below, each equal to 1 at $S_m$, with positive integer multipliers.*

- Seed: $S_0 = U(2k) - \Sigma[k..2k-1] + \Sigma[0..k-2]$; $S_1 = U(2k+1) - \Sigma[k..2k] + \Sigma[0..k-1]$; $S_2 = U(2k+2) - \Sigma[k+1..2k+1] + \Sigma[0..k-1]$; $S_3 = U(2k+3) - \Sigma[k+1..2k+2] + \Sigma[0..k]$; $S_4 = U(2k+4) - \Sigma[2k+1..2k+3] - \Sigma[k+1..2k-1] + \Sigma[0..k]$. These are $V_{m'}$ for $m' = 2k, \ldots, 2k + 4$.
- Recursion: $V_{m'}$ for $2k + 4 < m' \le m$, generated from $V_{2k+4} = S_4$ as in Section 4.2 ($i_0 = 2k + 4$, $\tau = \Sigma[0..k] - \Sigma[k+1..2k-1]$).
- Top block ($2k - 1$ relations): $D = U(k) - U(k-1)$; $P_h = U(h) - \Sigma[k..h-1] + \Sigma[2k-h..k-2]$ and $N_h = -\Sigma[k+1..h] + \Sigma[2k-h..k] - U(2k-1-h)$ for $k + 1 \le h \le 2k - 1$.

**Theorem 4.4 (Family II).** *Let $k \ge 2$ and $m \ge 2k + 1$. The same conclusion holds for the set $S_m$ of Family II, with the following m relations.*

- Seed: $T_0 = U(2k+1) - \Sigma[k..2k] + \Sigma[0..k-1]$; $T_1 = U(2k+2) - \Sigma[k+1..2k+1] + \Sigma[0..k-1]$; $T_2 = U(2k+3) - \Sigma[k+1..2k+2] + \Sigma[0..k]$; $T_3 = U(2k+4) - \Sigma[k+2..2k+3] + \Sigma[0..k]$; $T_4 = U(2k+5) - \Sigma[2k+2..2k+4] - \Sigma[k+2..2k] + \Sigma[0..k]$. These are $V_{m'}$ for $m' = 2k + 1, \ldots, 2k + 5$.
- Recursion: $V_{m'}$ for $2k + 5 < m' \le m$, generated from $T_4$ ($i_0 = 2k + 5$, $\tau = \Sigma[0..k] - \Sigma[k+2..2k]$).
- Top block ($2k$ relations): $D' = U(k) - U(k-1)$; $P'_h = U(h) - \Sigma[k..h-1] + \Sigma[2k+1-h..k-1]$ for $k + 1 \le h \le 2k$; $N'_h = -\Sigma[k+2..h] + \Sigma[2k+1-h..k] - U(2k-h)$ for $k + 2 \le h \le 2k$.

### *4.4 Proofs*

We give the proof of Theorem 4.3; the proof of Theorem 4.4 is identical, with the changes indicated at the end.

**Step 1: closed forms.** For Family I, $U(t) = (2/3)(4^{k-1} - 4^{k-1-t})$ for $0 \le t \le k - 1$, $U(k) = (2 \cdot 4^{k-1} + 1)/3$, and $U(t) = U(k) + (4^{t-k} - 1)/3$ for $k \le t \le 2k$. The rule values $b = k + 1, k + 1, k + 2, 3$ at $i = 2k + 1, \ldots, 2k + 4$ give $U(2k+1) = 2U(2k) - U(k-1)$, $U(2k+2) = 2U(2k+1) - U(k)$, $U(2k+3) = 2U(2k+2) - U(k)$, $U(2k+4) = 2U(2k+3) - U(2k)$.

**Step 2: tightness of the seed and of the block.** With Step 1, each of $S_0, \ldots, S_4, D, P_h, N_h$ is a finite combination of geometric sums in powers of 4. For example, with $q = h - k \in [1, k - 1]$,

$$P_h = U(h) - \Sigma_{t=k}^{h-1} U(t) + \Sigma_{t=2k-h}^{k-2} U(t) = (1 - q)U(k) + (4^q - 1)/3 - (4^q - 1)/9 + q/3 + (2/3)(q - 1)4^{k-1} - (2/9)(4^q - 4),$$

and since $(1 - q)U(k) + (2/3)(q - 1)4^{k-1} = (1 - q)/3$, the terms in $4^q$ cancel and $P_h = 1$. The other identities are computations of the same kind; all of them ($S_0, \ldots, S_4, D, P_h, N_h$, and $T_0, \ldots, T_4, D', P'_h, N'_h$) were also verified symbolically in k and h.

**Step 3: the recursion.** For $i \ge 2k + 4$ we have $b(i) = r(i + 2 - 2k)$ with $r(N) = \mathrm{round}(\sqrt{2N})$, and $b(2k + 4) = r(6) = 3$ is the block length of $S_4$. By Lemma 3.3, r increases in isolated unit steps for $N \ge 6$, the first increase being at $i = 2k + 5$. So Lemma 4.1 applies with $i_0 = 2k + 4$, and every $V_{m'}$ has value 1. Lemma 4.2 applies as well: $\tau = \Sigma[0..k] - \Sigma[k+1..2k-1]$ has coefficient sum $2 = b(2k+4) - 1$, its largest index $2k-1$ is smaller than the first new index $2k+5-1-4=2k$, and its partial sums from index 0 are positive. The seeds $S_0, \ldots, S_3$ have tails of $+1$ coefficients only, with sum $b-1$ ($b = k, k+1, k+1, k+2$), so they have the same triangular pattern.

**Step 4: the triangular part.** The relations $V_{m'}$, $2k \le m' \le m$, have leading gap positions $j = m - m' + 1 = 1, \ldots, m-2k+1$. As in Theorem 3.2, the equations $\Sigma$ y E = 1 on these columns give $y_j = 1 + \Sigma_{i<j}\, y_i \,|E^{(i)}{}_j| \ge 1$.

**Step 5: the block.** The block relations use only $U(0), \ldots, U(2k-1)$, so they vanish on the gap columns $1, \ldots, m-2k+1$; on the remaining $2k-1$ columns, which correspond to $d(1), \ldots, d(2k-1)$, a relation $\Sigma\, c_t\, U(t)$ has coefficient $e_i = \Sigma_{t\ge i}\, c_t$ at $d(i)$. A direct count gives: D has e = unit vector at $i = k$; $P_h$ has $e_h = 1$, $e_i = 0$ for $i > h$ and $e_i \le 0$ for $i < h$; $N_h$ has $e_{2k-h} = 1$, $e_i = 0$ for $i < 2k-h$ and $e_i \le 0$ for $i > 2k-h$. The right-hand sides on these columns are $t_i = 1 - \Sigma_{\text{triangular}}\, y\, E_i \ge 1$. Solve the columns in the order $1, 2k-1, 2, 2k-2, \ldots, k-1, k+1$, and finally k. At column $i \le k-1$ the only unknown is $N_{2k-i}$; at column $2k-i$ the only unknown is $P_{2k-i}$; at column k the only unknown is D. In each equation the new multiplier equals $t_i$ plus a nonnegative combination of multipliers already found, so all block multipliers are positive. This order also shows that the block matrix is triangular after a permutation, with determinant $\pm 1$.

**Step 6: conclusion.** The m relations are linearly independent (block triangular, with both diagonal blocks triangular), and $\Sigma\, y_j\, R_j = x_m$ with all $y_j > 0$. Their values at $S_m$ are 1, so $\Sigma\, y_j = \max S_m$ and $x_m - \max S_m = \Sigma\, y_j\, (R_j(x) - 1) \ge 0$ on the chamber, with equality only at $S_m$. For $2k \le m < 2k+4$ only $V_{2k}, \ldots, V_m$ are used and the same argument applies. ■

**Changes for Theorem 4.4.** $U(t) = (4^k - 4^{k-t})/3$ for $0 \le t \le k$, $U(k+1) = U(k) + 1$, $U(t) = U(k+1) + 2(4^{t-k-1} - 1)/3$ for $k+1 \le t \le 2k+1$, and $U(2k+2) = 2U(2k+1) - U(k)$, $U(2k+3) = 2U(2k+2) - U(k)$, $U(2k+4) = 2U(2k+3) - U(k+1)$, $U(2k+5) = 2U(2k+4) - U(2k+1)$. The recursion starts at $i_0 = 2k+5$ with $b(2k+5) = 3$; the first new index is $2k+1$. In the block, $D'$ leads at k, $P'_h$ at $h \in [k+1, 2k]$ (with $e_i = 0$ for $i > h$ and $e_i \le 0$ for $i < h$) and $N'_h$ at $2k+1-h \in [1, k-1]$ (with $e_i = 0$ below and $e_i \le 0$ above); the columns are solved in the order $1, 2k, 2, 2k-1, \ldots, k-1, k+2$, then $k+1$ and k.

**Computational confirmation.** The certificates generated from these formulas were checked in exact arithmetic (coefficients, tightness, triangular pattern, positivity, $\Sigma\, y_j = \max S_m$) for $2 \le k \le 12$ and all m up to $2k+40$, for both families. They agree with the chamber LPs solved directly with all $(3^m - 1)/2$ relations ($k = 2$, $12 \le m \le 16$; $k = 3$, $m = 16$).

Theorems 4.3 and 4.4 show that chamber optimality is not special to the Conway–Guy sets: it holds for the constructions that improve the Conway–Guy constant (Bohman's $0.22002 \cdot 2^n$). It is not automatic for every construction of the same type (see Section 7).

## 5. Exact values for $n \le 8$

**Theorem 5.1 (Theorem B).** *$f_R(n) = f_Z(n)$ for $1 \le n \le 8$: the values are 1, 2, 4, 7, 13, 24, 44, 84.*

The upper bounds are attained by the Conway–Guy sets. For the lower bounds we use a branch-and-bound over chambers.

### *5.1 The tree*

A node is a finite list of oriented relations $\rho_1, \ldots, \rho_k$. Its polyhedron consists of the base constraints $x_1 \ge 1$, $x_i - x_{i-1} \ge 1$, the valid cuts below, and $\rho_l(x) \ge 1$. We minimise $x_n$ over this polyhedron. If the minimum is at least $f_Z(n)$, or the polyhedron is empty, the node is a leaf. Otherwise the minimiser violates some relation ε ($|\varepsilon\cdot x| < 1$) and the node gets the two children obtained by adding $\varepsilon\cdot x \geq 1$ or $-\varepsilon\cdot x \geq 1$. Every admissible point satisfies one of these, so the leaves cover the admissible set.

**Valid cuts** (used for n = 7 and n = 8): $x_1 + \cdots + x_k \geq 2^k - 1$, because the $2^k$ subset sums of $\{x_1, \ldots, x_k\}$ lie in $[0, x_1 + \cdots + x_k]$ and are 1 apart; $x_k \geq f_R(k)$ for $k < n$, because $\{x_1, \ldots, x_k\}$ is admissible; and $x_n \geq H_n$ (Proposition 1.2). The cases $n \leq 6$ use no cuts: their polyhedra contain only the base constraints and the branching relations, so these trees depend on no earlier result. (For $n \leq 5$, Proposition 1.2 alone would already give the lower bound, since $H_n = f_Z(n)$ for $n \leq 5$.) The cases n = 7 and n = 8 use the cuts, with the values $f_R(k)$, $k < n$, proved by the trees for smaller n. Section 5.3 and Appendix A describe how this chain is checked, starting at n = 1.

### *5.2 Leaf certificates*

A leaf with system $Ax \geq b$, $x \geq 0$ is closed by a rational vector $y \geq 0$ with $A^T y \leq e_n$ componentwise and $b\cdot y \geq f_Z(n)$: then $x_n \geq y^T Ax \geq y\cdot b$ for every feasible $x \geq 0$. This also covers leaves with empty polyhedron. The vectors y were obtained from floating-point LP duals, converted to rationals, and accepted only after an exact check.

| n | nodes | leaves | cuts | result |
|---|---|---|---|---|
| 1 | 1 | 1 | none | $f_R(1) = 1$ |
| 2 | 1 | 1 | none | $f_R(2) = 2$ |
| 3 | 3 | 2 | none | $f_R(3) = 4$ |
| 4 | 13 | 7 | none | $f_R(4) = 7$ |
| 5 | 173 | 87 | none | $f_R(5) = 13$ |
| 6 | 4,683 | 2,342 | none | $f_R(6) = 24$ |
| 7 | 384,005 | 192,003 | §5.1 | $f_R(7) = 44$ |
| 8 | 112,785,751 | 56,392,876 | §5.1 | $f_R(8) = 84$ |

### *5.3 Verification*

Correctness rests on two properties: the leaves form a complete binary tree in which each internal node splits on a pair ±ε, and each leaf certificate satisfies $y \geq 0$, $A^T y \leq e_n$ and $b\cdot y \geq f_Z(n)$ exactly.

For $n \leq 8$ both properties were checked on the stored data by a program independent of the search code. It rebuilds every constraint system from the stored path and checks each certificate in exact integer arithmetic; coverage is checked by a streaming test that the leaves, listed in depth-first order, form a complete binary tree. The test is specified as follows.

**Nodes and complete trees.** A node is identified by its path $(\rho_1, \ldots, \rho_d)$, the oriented relations added on the way from the root, each a nonzero vector in $\{-1,0,1\}^n$. A tree is complete if every internal node p has exactly the two children $p + (\varepsilon)$ and $p + (-\varepsilon)$ for one such ε. The leaves of a subtree with root r of depth $d_0$ are stored in depth-first order: all leaves below one child of a node before all leaves below the other.

**Streaming test.** The checker reads the leaves once and keeps a stack of pending siblings $(l, \sigma)$. Such a pair means that the node $L[..l) + (\sigma)$ has still to be visited, where L is the leaf read last and $L[..l)$ consists of its first l relations. The first leaf $L_1$ pushes $(l, -L_1[l])$ for $d_0 \leq l < |L_1|$. Each following leaf $L'$ pops the top pair $(l, \sigma)$ and must have the same ancestral path as the previous leaf L, that is $L'[..l) = L[..l)$, and the opposite relation $L'[l] = \sigma$ at level l; it then pushes $(l', -L'[l'])$ for $l < l' < |L'|$. Every relation read is checked to be nonzero with entries in $\{-1,0,1\}$. At the end there must be at least one leaf and no pending sibling; a root that is itself a leaf must be the only leaf.

**Invariant.** After each leaf L, the stack contains, in increasing order of level, exactly the pairs $(l, -L[l])$ for which L lies below the child $L[..l) + (L[l])$ listed first and the other child has not been visited yet. The next leaf in depth-first order is therefore the first leaf below the deepest pending sibling, which is what the pop tests, and an empty stack at the end means that every internal node on the path of every leaf has both children. By induction on the number of leaves read, the test accepts exactly the depth-first leaf lists of complete trees with root r.

**Top tree.** Each tree is stored as a top part, with the leaves found above a fixed frontier depth and the frontier nodes $F_0$, $F_1$, …, and one file for each frontier node. The checker verifies, by a recursive test with the same condition on splits, that the top leaves and the frontier nodes are the leaves of a complete tree from the root; that the file of $F_i$ exists and records $F_i$ as its root and i as its index; that every leaf in it extends $F_i$; and that its leaves pass the streaming test from $F_i$. Grafting the subtrees onto the frontier nodes gives one complete tree from the root whose leaves are exactly the stored leaves, each with a checked certificate.

**Cover.** Every nonzero $\varepsilon \in \{-1,0,1\}^n$ gives $\varepsilon\cdot x = \Sigma_{i\in S}\, x_i - \Sigma_{i\in T}\, x_i$ for disjoint sets S and T, so an admissible x satisfies $\varepsilon\cdot x \geq 1$ or $-\varepsilon\cdot x \geq 1$ at every split. Following these choices from the root, x reaches a leaf all of whose constraints it satisfies, and the certificate of that leaf gives $x_n \geq f_Z(n)$.

**Tests and the chain.** A test script in the package makes corrupted copies of the data for n = 6 and n = 7 (a missing, duplicated or reordered leaf, a wrong sibling relation, a negative or halved multiplier, a row index out of range, a missing, empty or exchanged subtree file, a missing frontier node, a raised target, a wrong value $f_R(6)$ in a tree with cuts, and an object, added to intact data, whose loading would run code) and checks that the checker rejects each of them and accepts the unmodified data. A second script checks the whole chain n = 1, …, 8 with one command: it runs the checker on every tree and confirms that each value $f_R(k)$ used as a cut was certified by the tree for k (Appendix A).

An earlier n = 7 tree (386,775 nodes, 193,388 leaves, without the bandwidth cut) was also checked by a second, separately written checker, and the n = 7 data reported in the table were rechecked with it. These additional checks were performed with the use of Astra (OpenAI), an AI system different from the one used for the search program and the main checker.

**The case n = 8.** Proposition 1.2 gives only $78 \leq f_R(8) \leq 84$; the tree closes the remaining gap.

The tree was split at depth 10 into 978 subtrees and 12 early leaves, solved with HiGHS in about 9 hours on two cores, and stored in compressed form (about 1.1 GB).

The node counts in the table are totals of the whole tree: for n = 8, 112,784,750 nodes inside the 978 subtrees (including their roots), 989 internal nodes above the split, and 12 early leaves; in each case the total equals $2\cdot(\text{leaves}) - 1$, as it must for a complete binary tree.

The independent checker verified all 56,392,876 leaf certificates and the complete coverage of the tree, subtree by subtree, using a streaming test of the depth-first leaf order.

A first version of the checker accepted a subtree file containing no leaves; this defect, found in the additional checks with Astra (OpenAI), was corrected (an empty subtree now fails). The data for n = 1, …, 8 in the package were verified with the current checker, which differs from the corrected one only by the option for trees without cuts, by a restricted loader for the data files (Appendix A), and by its documentation; for n = 8 the data were rebuilt from the 37-part archive and checked again with it (25 October 2026).

The search programs and the checkers were written with AI assistance (see Acknowledgements). The verification does not rely on the AI system: correctness follows from the exact checks of the published data, which anyone can rerun.

### *5.4 A span inequality*

The following inequality is a valid cut of a different kind: it bounds the distance between two elements, not an element itself.

**Proposition 5.1.** *Let $s(k)$ be the least possible value of $y_k$ over all reals $0 = y_0 < y_1 < \cdots < y_k$ such that $|\varepsilon\cdot y| \geq 1$ for every nonzero $\varepsilon \in \{-1,0,1\}^{k+1}$ with $\varepsilon_0 + \cdots + \varepsilon_k = 0$. If $0 < x_1 < \cdots < x_n$ is admissible, then $x_j - x_i \geq s(j - i)$ for all $i < j$. Moreover $s(k) \leq f_R(k)$.*

*Proof.* Put $y_t = x_{i+t} - x_i$ for $0 \leq t \leq j - i$. If $\varepsilon_0 + \cdots + \varepsilon_k = 0$, then $\varepsilon\cdot y = \Sigma_t\, \varepsilon_t\, x_{i+t}$, which is a signed relation of x, so $|\varepsilon\cdot y| \geq 1$. For the second claim, if $a_1 < \cdots < a_k$ is admissible, then $y = (0, a_1, \ldots, a_k)$ is feasible, because $\varepsilon\cdot y = \Sigma_{t\geq 1}\, \varepsilon_t\, a_t$ is a signed relation of a or zero, and it is zero only if $\varepsilon_1 = \cdots = \varepsilon_k = 0$, which together with $\varepsilon_0 + \cdots + \varepsilon_k = 0$ forces $\varepsilon = 0$. ■

In words: the elements $x_i, \ldots, x_j$, shifted by $x_i$, have 1-separated sums over subsets of equal size. The values for $k \leq 7$, computed exactly (upper bounds by the explicit points, lower bounds by a branch-and-bound over these relations with exact rational certificates and an independent checker), are:

| k | 1 | 2 | 3 | 4 | 5 | 6 | 7 |
|---|---|---|---|---|---|---|---|
| s(k) | 1 | 2 | 4 | 7 | 13 | 22 | 40 |
| $f_R(k)$ | 1 | 2 | 4 | 7 | 13 | 24 | 44 |

Optimal points are (0, 1, 2, 4, 7, 13) for k = 5, (0, 1, 2, 10, 16, 19, 22) for k = 6 and (0, 1, 2, 10, 19, 34, 37, 40) for k = 7. Thus $s(k) = f_R(k)$ for $k \leq 5$, and $s(k) < f_R(k)$ for k = 6, 7. The inequality $x_j - x_i \geq s(j - i)$ is stronger than the base inequality $x_j - x_i \geq j - i$ for $j - i \geq 3$. In our search, however, adding it as a cut did not reduce the trees: the optimal points of the node LPs almost always satisfy it already (for n = 7 the number of nodes even grew by about 4%, because different LP solutions lead to different branching).

More generally, every valid linear inequality for admissible m-sets can be applied to every m-element subset of $\{x_1, \ldots, x_n\}$. The strongest such inequalities are the facets of the convex hull of all admissible m-sets; this hull is the convex hull of the vertices of all chambers plus the cone of sorted directions $0 \leq d_1 \leq \cdots \leq d_m$. We computed these facets for $m \leq 5$ (4, 15 and 109 facets for m = 3, 4, 5, each verified exactly on every chamber by a rational dual vector). Examples are $x_1 - x_2 + x_3 \geq 3$ and $2x_1 - 3x_2 + 2x_3 + x_4 \geq 10$. Applied to all subsets, they reduced the tree for n = 6 from 4,475 to 3,583 nodes and for n = 7 from 401,731 to 350,959 nodes (counted without certificates), while the LP time per node grew about tenfold. Inequalities that involve only a few elements therefore do not explain the depth of the trees; a substantial reduction requires a global inequality, as in Proposition 1.2.

## 6. Non-integral chamber optima

An optimum of a chamber LP is attained at a vertex $x = M^{-1}\mathbf{1}$, where M is a nonsingular matrix of n tight relations. If $|\det M| = 1$ the vertex is integral; if $|\det M| \geq 2$ it may or may not be. For the Conway–Guy chambers M is unit triangular in gap coordinates (Theorem 3.2), which explains why the real relaxation gives no improvement there. A strict inequality $f_R(n) < f_Z(n)$ requires a chamber with a non-integral optimal vertex that beats every integer chamber.

A local search over chambers near the Conway–Guy and Bohman optima ($8 \le n \le 12$) found non-integral chamber optima. Two examples, each verified exactly (admissibility, chamber membership, and a positive dual certificate equal to the maximum):

- $n = 11$: $x = (1455/2, 2173/2, 2553/2, 1394, 1439, 1471, 1482, 2969/2, 1489, 1507, 1526)$, chamber minimum 1526.
- $n = 12$: $x = (3581/2, 5373/2, 6333/2, 6839/2, 7069/2, 3607, 7285/2, 7289/2, 7329/2, 3689, 7393/2, 7413/2)$, chamber minimum 7413/2.

The same search also found an integral optimal vertex with $|\det M| = 2$ ($n = 12$, maximum 2092), so the determinant alone does not decide integrality.

All non-integral optima found are far above the best known values, and no chamber better than the known constructions was found. This search is heuristic.

For $n \le 5$ the optimum of every chamber is unique and is an integer point. The smallest non-integral examples occur at $n = 6$; Appendix B gives one with its matrix M and shows how $|\det M| = 2$ appears in the tight relations.

## 7. Open problems

1. Extend Theorems 4.3 and 4.4 to further constructions. The four sequences of Lunnon [L88] that improve on the Conway–Guy sequence are Families I and II with $k = 2$ and $k = 3$ [B98, §2], so Theorems 4.3 and 4.4 already cover them. Open are Lunnon's other generalizations of the Conway–Guy sequence [L88, §7], general CG-generators (Problem 2) and the sets of the 2026 disproof. Lemmas 4.1 and 4.2 reduce such a case to a finite seed and a top block only when their hypotheses hold: the recursion $d(i) = d(i - b(i)) + \cdots + d(i - 1)$ for $i > i_0$, a rule sequence b that increases in unit steps and never at two consecutive indices, and a tail $\tau$ with properties (ii) and (iii) of Lemma 4.2. The CG-generator of Problem 2 has consecutive increases of b, and its set is not chamber-optimal.
2. Bohman [B96, §3, Theorem 2] shows that any rule sequence b with $b(1) = 0$, $b(i + 1) \in \{b(i), b(i) + 1\}$ and $b(i - b(i)) < b(i)$ generates sets with distinct subset sums (CG-generators). Theorem 3.2 does not extend to all CG-generators: for Bohman's example $b(i) = i - 1$ ($i \le 6$), $b(i) = \mathrm{round}(\sqrt{2(i - 4)})$ ($i > 6$), the set with $m = 9$ is (88, 136, 160, 176, 184, 188, 190, 191, 192), but its chamber also contains the integer set (80, 124, 148, 160, 168, 172, 174, 175, 176), which has distinct subset sums and smaller maximum. Characterise the CG-generators whose sets are chamber optima.
3. Push the certified computation to $n = 9$ (where $f_Z(9) = 161$ and $H_9 = 148$). The tree grew by a factor of about 300 from $n = 7$ to $n = 8$. A random-dive estimate of the tree size [K75], with the cuts of Section 5.1 and $x_8 \ge f_R(8) = 84$, gives between $7 \cdot 10^9$ and $6 \cdot 10^{10}$ nodes for each target $t \in \{149, 150, 152, 155, 161\}$; the median depth of a leaf is 24–26 for all these targets, and the same estimator underestimated the $n = 8$ tree by a factor of about 7. So $n = 9$ needs stronger relaxations, not only more computing time.
4. *Conjecture.* $f_R(n) \ge H_n + 1$ for all $n \ge 6$. This holds for $n = 6, 7, 8$ ($24 > 23$, $44 > 43$, $84 > 78$), and for $n = 9$ it would give the first improvement $f_R(9) \ge 149$. The order of the subset sums alone does not suffice: for $3 \le n \le 9$ the order of the subset sums of the Conway–Guy set is a numbering of the cube of bandwidth exactly $H_n$ (computed), so Harper's bound is attained by a realizable order. For $n = 6$, all realizable orders found with bandwidth exactly 23 have chamber minimum at least 24, with equality only for the Conway–Guy chamber (Theorem 3.2).
5. Decide whether $f_R(n) = f_Z(n)$ for all n. Integer upper bounds (records) are not enough to show $f_R(n) < f_Z(n)$: a proof needs an integer lower bound as well, and determining the first such n requires excluding all smaller n.
6. Find a uniform certificate for a family of sets whose maxima are $o(u_n)$, and relate the certificate multipliers to the asymptotic constants of Conway–Guy and Bohman.

## Appendix A. Data and reproducibility

The programs and the certificate data are deposited at Zenodo: doi:10.5281/zenodo.23315221 (this DOI covers all versions and resolves to the latest) [No26].

The programs in the accompanying package, and the computations they carry out, were produced with AI assistance (Claude, Anthropic). They use Python 3 with NumPy, SciPy, HiGHS (via highspy), SymPy and python-flint. Exact checks use Python integers and rationals. The package was last run in full with Python 3.14.5, NumPy 2.5.4, SciPy 1.18.1, HiGHS 1.15.1 and SymPy 1.14.0.

The accompanying package contains:

- the branch-and-bound program, with an option for trees without cuts, and the independent checker (exit code 0 only for a certified result), which reads the data files with a restricted loader that admits only plain data (so checking downloaded data runs no code from them) and whose documentation contains the specification of Section 5.3;
- the complete certified data for $n = 1, \dots, 7$ (for $n \le 6$ without cuts), with the commands that produced them;
- a script that checks the whole chain with one command, `python3 verify_chain.py`: it runs the checker on the trees for $n = 1, \dots, 7$, and for $n = 8$ when those data are present, confirms that every value $f_R(k)$ used as a cut was certified earlier in the chain, and checks the upper bounds by the Conway–Guy sets; and a test script, `python3 test_checker.py`, with corrupted data that the checker must reject;
- for $n = 8$, the search and check logs and SHA-256 hashes of all 979 data files;
- the $n = 8$ data themselves (979 files, about 1.1 GB), distributed separately as a 37-part archive with SHA-256 hashes of every part and instructions for reassembly; they can also be regenerated with the package and checked in the same way;
- the scripts for all other computations (Theorem 3.2 for $n \le 150$, the chamber LPs, Theorems 4.3 and 4.4 including the symbolic tightness checks, Section 6, Appendix B (`appB_check.py`, `small_chambers_exact.py`, `small_chambers.py`), Section 5.4 (`span_certify.py`, `span_check.py`; `hull.py` and the facet lists `facets_3.txt`, `facets_4.txt`, `facets_5.txt`), the tree-size estimate of Section 7 (`estimate_tree.py`) and the CG-generator example), with one script that runs them all.

## Appendix B. A non-integral chamber optimum for $n = 6$

For $n \le 5$ there is no non-integral chamber optimum. We sampled random points $0 < x_1 < \cdots < x_n$ and collected their sign chambers. This found 2, 14 and 516 chambers for $n = 3, 4, 5$, which are all of them: these are the numbers of orders of the $2^n$ subset sums that arise from n positive reals [Mac99; OEIS A009997].

For each of these 532 chambers we certified the optimum in exact arithmetic, by n tight relations M with $x = M^{-1}\mathbf{1}$ in the chamber and multipliers $y = M^{-T}e_n \ge 0$. In every case the relations with $y_j > 0$ have rank n, so the optimum is unique, and it is an integer point. Hence an example with a $3 \times 3$, $4 \times 4$ or $5 \times 5$ matrix M does not exist.

For $n = 6$ the same sampling found 122,826 of the 124,187 chambers. Six of them have a non-integral optimal vertex, each with $|\det M| = 2$. The one with the smallest maximum is

$$p = (20, 45/2, 25, 61/2, 34, 38).$$

The six tight relations (the rows of M) and the multipliers y are:

| j | tight relation $R_j(x) = 1$ | $y_j$ |
|---|---|---|
| 1 | $x_1 - x_2 - x_4 + x_5$ | 9 |
| 2 | $x_1 + x_2 + x_4 - x_5 - x_6$ | 9 |
| 3 | $-x_1 - x_2 - x_3 + x_4 + x_6$ | 8 |
| 4 | $-x_1 + x_2 + x_3 - x_4 - x_5 + x_6$ | 6 |
| 5 | $-x_1 + x_3 + x_5 - x_6$ | 4 |
| 6 | $x_2 - x_3 - x_4 + x_5$ | 2 |

*Verification.* The point p is admissible: $|\varepsilon\cdot p| \geq 1$ for all 728 nonzero ε. Each $R_j$ equals 1 at p, so $R_j \geq 1$ is a constraint of the chamber C(p). The multipliers satisfy $\Sigma_j\, y_j\, R_j(x) = x_6$ identically and $\Sigma_j\, y_j = 38$.

Hence $x_6 - 38 = \Sigma_j\, y_j\, (R_j(x) - 1) \geq 0$ on C(p). Since all $y_j > 0$ and M is nonsingular, equality holds only at p. So p is the unique optimum of its chamber, and no integer point of the chamber attains the value 38, although 38 is an integer.

*How the determinant appears.* Here det M = −2, so $M^{-1}$ = adj(M)/det M and 2p is an integer vector. Put v = (0, 1, 0, 1, 0, 0). Every tight relation contains $x_2$ and $x_4$ together (rows 1, 2, 3, 4, 6) or neither of them (row 5), and

$$M\, v = (-2, 2, 0, 0, 0, 0) \equiv 0 \pmod 2.$$

So v is a nonzero kernel vector of M modulo 2, which is possible only because det M is even. Shifting $x_2$ and $x_4$ by 1/2 together changes every $R_j$ by an integer, and this is what allows the half-integral vertex p = z + v/2 with z = (20, 22, 25, 30, 34, 38).

The integer point z is not admissible: M z = (2, 0, 1, 1, 1, 1), and $R_2(z) = 0$ says 20 + 22 + 30 = 34 + 38. Rounding p down therefore creates two equal subset sums.

The chamber minimum 38 is far above $f_R(6) = 24$, in agreement with Theorem B. The example is checked in exact rational arithmetic by `appB_check.py`; the cases $n \leq 5$ by `small_chambers_exact.py`, and the sampling for n = 6 by `small_chambers.py` (Appendix A). `small_chambers_exact.py` exits with an error unless it finds exactly 2, 14 and 516 chambers for n = 3, 4, 5 and certifies every one of them with a unique, integral optimum and no failure. `small_chambers.py` exits with an error if a sampled chamber with a non-integral optimum is left uncertified, or if the documented run does not give 122,826 chambers with six non-integral optima.

## Acknowledgements

The computations, the programs and the main verifications reported in this paper were carried out with the assistance of an AI system (Claude, Anthropic). The checks are exact and are performed by the published programs, so they can be repeated independently of the AI system.

Additional checks were performed with the use of Astra (OpenAI): a separate checker for the n = 7 computation was written, the n = 7 certificates were re-checked, and a defect in an earlier version of our checker (acceptance of an empty subtree file) was found.